\documentclass{kybernetika}

\usepackage{mathtools}
\usepackage{array}
\usepackage{cite}
\usepackage[nameinlink,capitalise]{cleveref}
\hypersetup{hidelinks}

\newtheorem{theorem}{Theorem}[section]
\newtheorem{lemma}[theorem]{Lemma}
\newtheorem{proposition}[theorem]{Proposition}
\newtheorem{corollary}[theorem]{Corollary}
\newtheorem{remark}[theorem]{Remark}
\newtheorem{example}[theorem]{Example}
\newtheorem{definition}[theorem]{Definition}

\newcommand{\preceqT}{\mathrel{\trianglelefteq}}
\newcommand{\precT}{\mathrel{\triangleleft}}
\newcommand{\npreceqT}{\mathrel{\not\trianglelefteq}}
\newcommand{\Ie}{I_e}
\newcommand{\cl}{\operatorname{cl}}
\newcommand{\intr}{\operatorname{int}}
\newcommand{\rtr}{\mathrm{rtr}}
\newcommand{\ltr}{\mathrm{ltr}}
\newcommand{\mtr}{\mathrm{mtr}}
\newcommand{\tr}{\mathrm{tr}}
\begin{document}
\pagestyle{myheadings}

\title{\texorpdfstring{Uninorm constructions on bounded trellises
via closure and interior operators}%
{Uninorm constructions on bounded trellises via closure and interior operators}}

\author{  Zhenyu Xiu}

\contact{Zhenyu}{Xiu}{Corresponding author. School of Mathematical Science, Guangxi Minzu University, Nanning 530006, Guangxi, China}{xyz198202@163.com}


\maketitle

\begin{abstract}

In this paper,  we establish necessary and sufficient conditions under which six closure- or
interior-based operations on bounded trellises are uninorms.  These operations
incorporate triangular norms, triangular conorms, or both, and form three
order-dual pairs.  Since the underlying pseudo-order need not be transitive
and meet and join need not be associative, standard lattice arguments do not
apply directly.  Under a fixed ambient transitivity assumption, the first pair
is characterized by image, monotonicity, and compatibility conditions for the
associated maps on elements incomparable with the neutral element; an
additional comparison condition used in earlier sufficient constructions is
shown to be redundant.  Complete characterizations are also obtained for the
variants with an absorbing endpoint, including both nondegenerate and
degenerate interval cases.  For the pair involving both a triangular norm and
a triangular conorm, we derive exact comparison conditions requiring the
triangular norm to take the value~$0$, or dually the triangular conorm to take
the value~$1$, in the relevant cases.  These conditions reduce to the usual
incomparability assumptions in the absence of nontrivial zero- or
one-divisors.  On bounded lattices, the closure and interior conditions
simplify to explicit identities or image restrictions.  Finite proper-trellis
examples clarify the role and limitations of the ambient transitivity
assumptions.

\end{abstract}



\section{INTRODUCTION}\label{sec:introduction}

Uninorms were introduced by Yager and Rybalov as associative and commutative
aggregation operations whose neutral element need not be an endpoint
\cite{YagerRybalov1996}.  Their structure on the unit interval has been studied
from several viewpoints, including representability and ordinal sums
\cite{Fodor1997,Drygas2016,MesiarovaZemanekova2017}.  More broadly,
lattice-valued fuzzy-set frameworks have long motivated topological and
logical investigations \cite{Goguen1967,Goguen1973}.

The extension of uninorms to bounded lattices led to a variety of construction
and characterization results \cite{Bodjanova2014,KaracalMesiar2015,Cayli2016,
Cayli2019,Dan2019,ZhaoWu2021,Sun2022,Zhang2021,Li2021,Qin2023,
MesiarovaZemanekova2024}.  Of particular relevance here are constructions and
characterizations based on closure and interior operators
\cite{OuyangZhang2020,CayliKybernetika2023,QiXiu2025,Xiu2025}.  These operators
are useful when defining the values of an operation on elements incomparable
with its neutral element.

A bounded trellis replaces the transitive order of a bounded lattice by a
reflexive and antisymmetric pseudo-order while retaining a meet and a join for
each pair of elements \cite{fried1970tournaments,Skala1971,
gladstien1973characterization}.  In this setting, meet and join need not be
associative, and consecutive comparisons cannot be concatenated without a
suitable transitivity hypothesis.  Zedam and De Baets \cite{ZedamDeBaets2023}
studied triangular norms on bounded trellises.  Kong and Zhao
\cite{KongZhao2024} subsequently introduced uninorms on bounded trellises and
gave constructions using triangular norms, triangular conorms, and closure or
interior operators.  Further developments appear in
\cite{LuoPangYao2025,WangDanDeBaets2026,SongZhao2026,XiuZheng2025,
ZhaoLiu2025,Mehenni2023}.

This paper examines six closure- or interior-based formulas, arranged in three
order-dual pairs.  Its aim is not to introduce another unrestricted
construction scheme, but to determine when these particular formulas define
uninorms.  Kong and Zhao \cite{KongZhao2024} gave several sufficient
constructions on bounded trellises, distinguished by transitivity properties
of the neutral element.  For the first pair considered here, (C1) is retained
as an ambient hypothesis, while (C2) and (C4) are shown to be necessary and
sufficient; (C3) then follows from (C1) and (C2).  The dual statement holds for
the interior-based formula.  The variants with an absorbing endpoint are
treated separately according to whether the relevant strict interval is empty.

For the pair involving both a triangular norm and a triangular conorm, the
usual incomparability assumption remains sufficient, but the exact requirement
is the implication in \eqref{eq:zero-comparison-lower} or its order dual
\eqref{eq:one-comparison-upper}.  These conditions coincide with the
incomparability assumptions when there are no nontrivial zero- or one-divisors,
and a counterexample shows that the converse fails without this restriction.
On bounded lattices, the closure and interior parts reduce to explicit
identities or image restrictions.  The related lattice constructions in
\cite{OuyangZhang2020,CayliKybernetika2023,QiXiu2025,Xiu2025} use formulas or
input operations that are not identical in general.  The lattice results below
therefore compare the criteria for the present formulas rather than claim that
all of those constructions are special cases.

\Cref{sec:preliminaries} introduces the notation used below.
\Cref{sec:compatibility} collects the hypotheses on the closure and interior
data.  The six operations are studied in \Cref{sec:constructions}, and their
bounded-lattice forms are obtained in \Cref{sec:reduction}.

\section{PRELIMINARIES}\label{sec:preliminaries}

\begin{definition}\cite{fried1970tournaments}
A \emph{pseudo-order} on a nonempty set $X$ is a reflexive and antisymmetric binary relation $\preceqT$.  A pseudo-ordered set is denoted by $(X,\preceqT)$.  We write $x\precT y$ when $x\preceqT y$ and $x\ne y$, and $x\parallel y$ when neither $x\preceqT y$ nor $y\preceqT x$ holds.
A pseudo-ordered set is called \emph{bounded} if it has a least element $0$ and a greatest element $1$.
\end{definition}

\begin{definition}\cite{gladstien1973characterization}
A \emph{trellis} is a pseudo-ordered set $(X,\preceqT)$ in which every pair $x,y$ has a greatest lower bound $x\wedge y$ and a least upper bound $x\vee y$.  It is \emph{bounded} if it has a least element $0$ and a greatest element $1$.  A bounded trellis is written $\mathbb T=(X,\preceqT,\wedge,\vee,0,1)$.
\end{definition}

The relation of a trellis is not assumed to be transitive.  Accordingly, all
iterated meets and joins in the trellis arguments below are explicitly
parenthesized.

The order dual of a bounded trellis
$\mathbb T=(X,\preceqT,\wedge,\vee,0,1)$ is
$\mathbb T^{\partial}=(X,\preceqT^{\partial},\vee,\wedge,1,0)$, where
$x\preceqT^{\partial}y$ if and only if $y\preceqT x$
\cite{DaveyPriestley2002}.  It is again a bounded trellis.  Under order
duality, meet and join, the least and greatest elements, right- and
left-transitivity, closure and interior operators, and triangular norms and
triangular conorms are interchanged, whereas middle-transitivity is preserved.

Following \cite{KongZhao2024}, for $a\preceqT b$ put
$[a,b]=\{x\in X:a\preceqT x\preceqT b\}$ and
$[a,b)=\{x\in[a,b]:x\ne b\}$, and define $(a,b]$ and $(a,b)$
analogously. For a fixed $e\in X\setminus\{0,1\}$, the incomparable set is
$\Ie=\{x\in X:x\parallel e\}$.
Thus $X$ is the disjoint union of $[0,e)$, $\{e\}$, $(e,1]$, and $\Ie$.

\begin{definition}\cite{gladstien1973characterization}
A pseudo-ordered set $(X,\preceqT)$ is a $\wedge$-semi-trellis
(respectively, a $\vee$-semi-trellis) if every pair of elements has a
greatest lower bound (respectively, a least upper bound). A subset $A$ of a
trellis is a \emph{sub-trellis} if $x\wedge y,x\vee y\in A$ for all
$x,y\in A$, with the operations and pseudo-order restricted to $A$.
\end{definition}

\begin{definition}\cite{KongZhao2024}
An element $a\in X$ is right-transitive if $a\preceqT x\preceqT y$ implies $a\preceqT y$; it is left-transitive if $x\preceqT y\preceqT a$ implies $x\preceqT a$; and it is middle-transitive if $x\preceqT a\preceqT y$ implies $x\preceqT y$.  The corresponding sets are denoted by $X^{\rtr}$, $X^{\ltr}$, and $X^{\mtr}$.  Their intersection is denoted by $X^{\tr}$.
\end{definition}

\begin{definition}\cite{ZedamDeBaets2023}
A binary operation $F:X^2\to X$ is \emph{increasing} if
\[
  x_1\preceqT x_2\ \text{and}\ y_1\preceqT y_2
  \quad\Longrightarrow\quad
  F(x_1,y_1)\preceqT F(x_2,y_2).
\]
More generally, for $D\subseteq X^2$, a map $F:D\to X$ is increasing if the same implication holds whenever $(x_1,y_1),(x_2,y_2)\in D$.  A unary map $f:A\to X$, where $A\subseteq X$, is increasing if $x\preceqT y$ in $A$ implies $f(x)\preceqT f(y)$.
A triangular norm on a bounded pseudo-ordered set is an increasing, associative, and commutative operation with neutral element $1$.  A triangular conorm is defined dually, with neutral element $0$.  A \emph{uninorm} on $\mathbb T$ is an increasing, associative, and commutative operation $U:X^2\to X$ having a neutral element $e\in X$.
\end{definition}

The four-variable definition of increasingness is essential.  In a nontransitive pseudo-order, separate inequalities
$F(x_1,y_1)\preceqT F(x_2,y_1)$ and
$F(x_2,y_1)\preceqT F(x_2,y_2)$ do not by themselves imply the required joint comparison.

\begin{definition}\cite{OuyangZhang2020}
A map $\cl:X\to X$ is a \emph{closure operator} if it is extensive, increasing, and idempotent:
\[
 x\preceqT\cl(x),\qquad
 x\preceqT y\Rightarrow \cl(x)\preceqT\cl(y),\qquad
 \cl(\cl(x))=\cl(x).
\]
An \emph{interior operator} $\intr:X\to X$ is increasing and idempotent and satisfies $\intr(x)\preceqT x$.
\end{definition}

We first recall a necessary condition for the neutral element of a uninorm.

\begin{lemma}\cite{KongZhao2024}\label{lem:neutral-middle}
If a uninorm on a bounded trellis has neutral element $e$, then $e\in X^{\mtr}$.
\end{lemma}


\begin{proposition}\cite{ZedamDeBaets2023}
Let $\mathbb T$ be a bounded trellis, and let $T$ be a triangular norm on
the sub-trellis $[0,e]$ with neutral element $e$. Then $T(0,x)=0$ for every
$x\in[0,e]$.
\end{proposition}


\section{CONDITIONS ON THE CLOSURE AND INTERIOR DATA}\label{sec:compatibility}

Let $e\in X\setminus\{0,1\}$.  For a closure operator $\cl$, put
\begin{equation}\label{eq:ce}
  c_e(x)=\cl(x)\wedge(x\vee e),\qquad x\in\Ie,
\end{equation}
and
\begin{equation}\label{eq:me}
  m_e(x,y)=(x\wedge y)\wedge e, \qquad x,y\in\Ie.
\end{equation}

For the closure-based constructions, we use the following conditions.
\begin{enumerate}
\item[(C1)] $e\in X^{\mtr}\cap X^{\rtr}$ and $[0,e)\cup\Ie\subseteq X^{\rtr}$.
\item[(C2)] $c_e(\Ie)\subseteq\Ie$, the map $c_e$ is increasing on $\Ie$, and $c_e(c_e(x))=c_e(x)$ for every $x\in\Ie$.
\item[(C3)] If $x\in\Ie$, $y\in(e,1]$, and $x\preceqT y$, then $c_e(x)\preceqT y$.
\item[(C4)] $m_e:\Ie^2\to X$ is increasing in the joint sense, that is,
\[
x_1\preceqT x_2,\ y_1\preceqT y_2
\quad\Longrightarrow\quad
m_e(x_1,y_1)\preceqT m_e(x_2,y_2)
\qquad(x_i,y_i\in\Ie),
\]
and $m_e(c_e(x),y)=m_e(x,y)$ for all $x,y\in\Ie$.
\end{enumerate}
The corresponding interior-side conditions are stated explicitly.  For an interior operator $\intr$, put
\begin{equation}\label{eq:de}
 d_e(x)=\intr(x)\vee(x\wedge e),\qquad x\in\Ie,
\end{equation}
and
\begin{equation}\label{eq:ne}
  n_e(x,y)=(x\vee y)\vee e,\qquad x,y\in\Ie.
\end{equation}

\begin{enumerate}
\item[(I1)] $e\in X^{\mtr}\cap X^{\ltr}$ and $(e,1]\cup\Ie\subseteq X^{\ltr}$;
\item[(I2)] $d_e(\Ie)\subseteq\Ie$, the map $d_e$ is increasing on $\Ie$, and $d_e(d_e(x))=d_e(x)$ for every $x\in\Ie$;
\item[(I3)] if $y\in[0,e)$, $x\in\Ie$, and $y\preceqT x$, then $y\preceqT d_e(x)$.
\item[(I4)] $n_e:\Ie^2\to X$ is increasing in the joint sense, namely,
\[
x_1\preceqT x_2,\ y_1\preceqT y_2
\quad\Longrightarrow\quad
n_e(x_1,y_1)\preceqT n_e(x_2,y_2)
\qquad(x_i,y_i\in\Ie),
\]
and $n_e(d_e(x),y)=n_e(x,y)$ for all $x,y\in\Ie$.
\end{enumerate}

\begin{lemma}\label{lem:images}
For $x\in\Ie$, one has $x\preceqT c_e(x)$ and $d_e(x)\preceqT x$.
If \textnormal{(C1)--(C4)} hold, then $m_e(x,y)\in[0,e)$ for all
$x,y\in\Ie$.  Dually, \textnormal{(I1)--(I4)} imply
$n_e(x,y)\in(e,1]$.
\end{lemma}

\begin{Proof}
The element $x$ is a lower bound of $\cl(x)$ and $x\vee e$, so the defining property of the meet gives $x\preceqT c_e(x)$.  The dual argument gives $d_e(x)\preceqT x$.

Clearly $m_e(x,y)\preceqT e$.  If $m_e(x,y)=e$, then $e\preceqT x\wedge y\preceqT x$.  Right-transitivity of $e$ would give $e\preceqT x$, contrary to $x\in\Ie$. Here we use the hypothesis (C1) which says $e\in X^{\rtr}$. Thus $m_e(x,y)\in[0,e)$.  The last assertion is dual.
\end{Proof}

\begin{lemma}\label{lem:transitions}
Assume (C1) and let $x\preceqT y$.  Then
\begin{enumerate}
\item $x\in(e,1]$ implies $y\in(e,1]$;
\item $x=e$ implies $y\in\{e\}\cup(e,1]$;
\item $x\in\Ie$ implies $y\in\Ie\cup(e,1]$;
\item $y\in[0,e)$ implies $x\in[0,e)$.
\end{enumerate}
The order-dual statements hold under (I1).
\end{lemma}

\begin{Proof}
For the first two assertions use right-transitivity of $e$.  If $x\in\Ie$ and $y\preceqT e$, then $x\preceqT y\preceqT e$ and right-transitivity of $x$ contradicts $x\parallel e$, proving the third assertion.  For the fourth assertion, if $x\in(e,1]$, then item~(1) would force $y\in(e,1]$, contrary to $y\in[0,e)$; if $x=e$ or $x\in\Ie$, items~(2) and~(3), respectively, give the same contradiction.  Hence $x\in[0,e)$.  The uses of antisymmetry are confined to the endpoint identifications in the interval definitions.
\end{Proof}

\section{CONSTRUCTION METHODS FOR UNINORMS ON BOUNDED TRELLISES}\label{sec:constructions}

\subsection{A meet-based construction and its dual}

\begin{theorem}\label{thm:U1}\label{prop:U1-relative-characterization}
Let $\mathbb{T}$ be a bounded trellis, and let $e\in X\setminus\{0,1\}$.
Let $\cl$ be a closure operator, and let $S$ be a triangular conorm on the
sub-trellis $[e,1]$ with neutral element $e$.  Assume \textnormal{(C1)} and
define $U_1:X^2\to X$ by
\begin{equation}\label{eq:U1}
U_1(x,y)=
\begin{cases}
S(x,y), & x,y\in[e,1],\\
x, & (x,y)\in\bigl([0,e)\times[e,1]\bigr)
      \cup\bigl(\Ie\times\{e\}\bigr),\\
y, & (x,y)\in\bigl([e,1]\times[0,e)\bigr)
      \cup\bigl(\{e\}\times\Ie\bigr),\\
\cl(x)\wedge(x\vee e), & x\in\Ie,\ y\in(e,1],\\
\cl(y)\wedge(y\vee e), & y\in\Ie,\ x\in(e,1],\\
(x\wedge y)\wedge e, & x,y\in\Ie,\\
0, & \text{otherwise}.
\end{cases}
\end{equation}
Then the following statements are equivalent:
\begin{enumerate}
\item $U_1$ is a uninorm on $\mathbb T$ with neutral element $e$;
\item conditions \textnormal{(C2)} and \textnormal{(C4)} hold.
\end{enumerate}
Under these conditions, \textnormal{(C3)} follows from \textnormal{(C1)} and
\textnormal{(C2)}, and $0$ is an absorbing element of $U_1$.
\end{theorem}

\begin{Proof}
Assume first that \textnormal{(C2)} and \textnormal{(C4)} hold.  If
$x\in\Ie$, $y\in(e,1]$, and $x\preceqT y$, then $y$ is an upper bound of
$x$ and $e$.  Hence $x\vee e\preceqT y$.  Since
$c_e(x)\preceqT x\vee e$ and
$c_e(x)\in\Ie\subseteq X^{\rtr}$, right-transitivity gives
$c_e(x)\preceqT y$.  Thus \textnormal{(C3)} follows, and
\textnormal{(C1)--(C4)} are available in the sufficiency proof below.

The displayed clauses assign exactly one value to every pair.  The formula is
symmetric, and direct substitution gives $U_1(e,x)=x$ for every $x\in X$.

We first verify joint increasingness.  Take $x\preceqT x'$ and
$y\preceqT y'$.  Simultaneously interchanging
$(x,x')$ with $(y,y')$ leaves the required comparison unchanged, so
the roles of $x$ and $y$ may be interchanged whenever commutativity is used
below.
By \Cref{lem:transitions}, a point of $(e,1]$ can move only to
$(e,1]$, a point equal to $e$ can move only to $e$ or $(e,1]$, and a point
of $\Ie$ can move only to $\Ie$ or $(e,1]$.  We now inspect every possible
case for $(x,y)$.

If $(x,y)$ lies in $[0,e)^2$, $[0,e)\times\Ie$, or
$\Ie\times[0,e)$, then $U_1(x,y)=0$, so the required comparison follows
from the minimality of $0$.

Suppose $x\in[0,e)$ and $y\in[e,1]$.  Then $y'\in[e,1]$ and
$U_1(x,y)=x$.  If $x'\in[0,e)$, then
$U_1(x',y')=x'$ and $x\preceqT x'$.  If $x'\in[e,1]$, increasingness of
$S$ gives $x'=S(x',e)\preceqT S(x',y')=U_1(x',y')$.
Since $x\in[0,e)\subseteq X^{\rtr}$ by (C1) and $x\preceqT x'$, it follows that
$x\preceqT U_1(x',y')$.  If $x'\in\Ie$, then
\[
 U_1(x',y')=
 \begin{cases}
 x',&y'=e,\\
 c_e(x'),&y'\in(e,1],
 \end{cases}
\]
and $x'\preceqT c_e(x')$ by \Cref{lem:images}.  Right-transitivity of
$x$ again gives $x\preceqT U_1(x',y')$.

If $x,y\in[e,1]$, then $x',y'\in[e,1]$, and the desired inequality is
precisely the joint increasingness of $S$.

Suppose next that $x=e$ and $y\in\Ie$.  Here
$x'\in\{e\}\cup(e,1]$ and $y'\in\Ie\cup(e,1]$.  If $x'=e$, then
$U_1(x',y')=y'$, so $y\preceqT U_1(x',y')$.  If
$x'\in(e,1]$ and $y'\in\Ie$, then
$y\preceqT y'\preceqT c_e(y')=U_1(x',y')$,
and right-transitivity of $y$ yields the conclusion.  Finally, if
$x',y'\in(e,1]$, increasingness of $S$ gives
$y'=S(e,y')\preceqT S(x',y')=U_1(x',y')$.
Combining this with $y\preceqT y'$ by right-transitivity of $y$ proves the
remaining comparison.

Suppose $x\in(e,1]$ and $y\in\Ie$.  Then $x'\in(e,1]$ and
$y'\in\Ie\cup(e,1]$.  If $y'\in\Ie$, (C2) gives
$U_1(x,y)=c_e(y)\preceqT c_e(y')=U_1(x',y')$.  If $y'\in(e,1]$, then
(C3) and increasingness of $S$ give
$c_e(y)\preceqT y'=S(e,y')\preceqT S(x',y')$.
Because $c_e(y)\in\Ie\subseteq X^{\rtr}$, the two displayed comparisons
yield $U_1(x,y)\preceqT U_1(x',y')$ without invoking transitivity of the
ambient relation.

It remains to take $x,y\in\Ie$.  Put $r=m_e(x,y)$.  By
\Cref{lem:images}, $r\in[0,e)\subseteq X^{\rtr}$.  If
$x',y'\in\Ie$, (C4) directly gives
$r\preceqT m_e(x',y')$.  If, say, $x'\in\Ie$ and
$y'\in(e,1]$, then $U_1(x',y')=c_e(x')$.  From the defining properties of
the meets and repeated applications of right-transitivity of $r$, we obtain
\[
 r\preceqT x\wedge y\preceqT x\preceqT x'
 \preceqT c_e(x'),
\]
where each adjacent comparison is combined with the preceding one using
that the fixed initial element $r$ is right-transitive.  The case
$x'\in(e,1]$, $y'\in\Ie$ is symmetric.  If
$x',y'\in(e,1]$, the same argument first gives $r\preceqT x'$, while
$x'=S(x',e)\preceqT S(x',y')$.
One more application of right-transitivity of $r$ yields
$r\preceqT S(x',y')$.  The listed cases exhaust all transitions allowed by
\Cref{lem:transitions}; hence $U_1$ is jointly increasing.

It remains to verify associativity.  The element $e$ is neutral, and
inspection of \eqref{eq:U1} shows that $0$
is absorbing.  Thus triples containing $e$ or $0$ are immediate.  For all
remaining triples, the following list checks every possible choice of the
inner pair, not merely one ordering for each placement of the three arguments
among $[0,e)$, $(e,1]$, and $\Ie$.
\begin{enumerate}
\item If all three elements lie in $[0,e)$, every inner product is $0$, so
every bracketing is $0$.

\item If all three lie in $(e,1]$, the assertion is the associativity of
$S$.

\item If all three lie in $\Ie$, every inner product belongs to $[0,e)$ by
\Cref{lem:images}, and its product with the remaining incomparable element
is $0$.  Hence all three choices of the inner pair give $0$.

\item For two elements of $[0,e)$ and one element of $(e,1]$, pairing the
two lower elements first gives $0$, whereas pairing a lower and an upper
element first returns that lower element, whose product with the other
lower element is $0$.

\item For two elements of $[0,e)$ and one element of $\Ie$, every choice of
the inner pair gives $0$ either at the first step or at the second step.

\item Let $a\in[0,e)$ and $b,c\in(e,1]$.  Pairing $b,c$ first gives
$U_1(a,S(b,c))=a$.  Pairing $a$ with either upper element first also gives
$a$, and pairing the result with the remaining upper element again gives
$a$.

\item For one element of $[0,e)$ and two elements of $\Ie$, pairing the two
incomparable elements first gives an element of $[0,e)$ and hence $0$ after
the second multiplication; pairing the lower element with either
incomparable element first gives $0$ immediately.  Thus all bracketings are
$0$.

\item Let $a\in(e,1]$ and $b,c\in\Ie$.  Pairing $b,c$ first gives
$U_1(a,m_e(b,c))=m_e(b,c)$.  The other two choices of inner pair give,
respectively, $m_e(c_e(b),c)=m_e(b,c)$ and
$m_e(c_e(c),b)=m_e(c,b)=m_e(b,c)$,
where (C4) and the symmetry of $m_e$ are used.

\item Let $a\in\Ie$ and $b,c\in(e,1]$.  Pairing $b,c$ first gives
$c_e(a)$.  Pairing $a$ with either upper element first gives $c_e(a)$, and
the second multiplication gives $c_e(c_e(a))=c_e(a)$ by (C2).

\item Let $a\in[0,e)$, $b\in(e,1]$, and $c\in\Ie$.  Pairing $a,b$
first gives $a$ and then $0$; pairing $a,c$ first gives $0$; and pairing
$b,c$ first gives $c_e(c)$, whose product with $a$ is $0$.  Hence all three
bracketings are $0$.
\end{enumerate}
These are the ten possible placements of three arguments in $[0,e)$,
$(e,1]$, and $\Ie$, and every possible inner pair has been covered.
Therefore $U_1$ is associative and hence is a uninorm
with neutral element $e$.

Conversely, suppose that $U_1$ is a uninorm.  We first record two facts that
do not use \textnormal{(C2)--(C4)}.  For $x\in\Ie$,
\Cref{lem:images} gives $x\preceqT c_e(x)$.  Moreover,
for $x,y\in\Ie$ one has $m_e(x,y)\preceqT e$.  Equality
$m_e(x,y)=e$ would give $e\preceqT x\wedge y\preceqT x$,
and right-transitivity of $e$ would imply $e\preceqT x$, contrary to
$x\parallel e$.  Thus
\begin{equation}\label{eq:me-lower-under-C1}
 m_e(x,y)\in[0,e)\qquad(x,y\in\Ie).
\end{equation}

Fix $x\in\Ie$.  If
$c_e(x)\preceqT e$, then right-transitivity of $x$, applied to
$x\preceqT c_e(x)\preceqT e$, would give $x\preceqT e$, a contradiction.
Consequently, the decomposition
$X=[0,e)\cup\{e\}\cup(e,1]\cup\Ie$ shows that $c_e(x)$ belongs either to
$\Ie$ or to $(e,1]$.

Assume temporarily that $c_e(x)\in(e,1]$.  Associativity at the triple
$(x,1,1)$ gives
\[
 U_1\bigl(U_1(x,1),1\bigr)
 =U_1\bigl(x,U_1(1,1)\bigr).
\]
Since a triangular conorm on $[e,1]$ satisfies $S(t,1)=1$, the two sides
are $1$ and $c_e(x)$, respectively.  Indeed, increasingness and neutrality
give $1=S(e,1)\preceqT S(t,1)\preceqT1$,
so antisymmetry yields $S(t,1)=1$.  Hence $c_e(x)=1$.  Associativity at
$(x,x,1)$ would then give
\[
 m_e(x,x)
 =U_1\bigl(U_1(x,x),1\bigr)
 =U_1\bigl(x,U_1(x,1)\bigr)=1,
\]
contradicting \eqref{eq:me-lower-under-C1}.  Therefore
\begin{equation}\label{eq:ce-image-necessary}
 c_e(x)\in\Ie\qquad(x\in\Ie).
\end{equation}

If $x,y\in\Ie$ and $x\preceqT y$, increasingness of $U_1$, applied
directly to $(x,1)\preceqT(y,1)$, yields
$c_e(x)=U_1(x,1)\preceqT U_1(y,1)=c_e(y)$.
Thus $c_e$ is increasing on $\Ie$.  Using
\eqref{eq:ce-image-necessary}, associativity at $(x,1,1)$ now reads
\[
 c_e(c_e(x))
 =U_1\bigl(U_1(x,1),1\bigr)
 =U_1\bigl(x,U_1(1,1)\bigr)=c_e(x).
\]
This proves (C2).

On $\Ie^2$ the operation $U_1$ agrees with $m_e$.  Hence the joint
increasingness of $U_1$ implies the joint increasingness of $m_e$.  Finally,
for $x,y\in\Ie$, associativity at $(1,x,y)$ and
\eqref{eq:me-lower-under-C1} give
\[
\begin{aligned}
 m_e(c_e(x),y)
 &=U_1\bigl(U_1(1,x),y\bigr)\\
 &=U_1\bigl(1,U_1(x,y)\bigr)\\
 &=m_e(x,y).
\end{aligned}
\]
Thus (C4) also holds.
\end{Proof}

Kong and Zhao \cite{KongZhao2024} introduced uninorms on bounded trellises and
gave four construction methods based on closure or interior operators together
with suitable transitivity assumptions.  \Cref{thm:U1} has a more specific
purpose: after fixing the formula \eqref{eq:U1} and the ambient hypothesis
\textnormal{(C1)}, it identifies \textnormal{(C2)} and \textnormal{(C4)} as
the exact remaining conditions and shows that \textnormal{(C3)} is redundant.
In particular, the necessity proof derives the image, increasingness,
idempotence, and compatibility requirements from the uninorm axioms instead
of imposing them only as hypotheses.  This is the precise sense in which the
result is sharper for the fixed formula \eqref{eq:U1}; the necessity is part
of the conclusion rather than a restatement of the assumptions.
It is not a characterization of all closure-based uninorms on bounded
trellises, and the same distinction applies to \Cref{thm:U2} by order duality.
The following example shows both why \textnormal{(C1)} is used in
\Cref{thm:U1} and why it need not hold for every individual operation defined
by \eqref{eq:U1}.

\begin{example}
\label{ex:lattice-equality-insufficient}
\label{ex:C1-not-necessary}
Let $Z=\{0,a,b,e,1\}$.  Besides reflexivity and the boundary comparisons
$0\preceqT z\preceqT1$, impose only $a\preceqT b$ and $b\preceqT e$.
Thus $a\parallel e$.  Comparable pairs have their usual meet and join, and
$a\wedge e=0$, $a\vee e=1$.  This gives a bounded trellis, which is proper
because $a\preceqT b\preceqT e$ but $a\npreceqT e$.

Let $\cl$ be the identity map.  Then $\Ie=\{a\}$,
$[0,e)=\{0,b\}$, $(e,1]=\{1\}$, $c_e(a)=a$, and $m_e(a,a)=0$.
In particular, the lattice-style equality
$a\wedge e=\cl(a)\wedge e$ holds.  On $[e,1]=\{e,1\}$ take the unique
triangular conorm.  Formula \eqref{eq:U1}, in the order $0,a,b,e,1$, becomes
\[
\begin{array}{c|ccccc}
 U_1&0&a&b&e&1\\ \hline
 0&0&0&0&0&0\\
 a&0&0&0&a&a\\
 b&0&0&0&b&b\\
 e&0&a&b&e&1\\
 1&0&a&b&1&1
\end{array}.
\]
The table is symmetric and has neutral element $e$.  Associativity follows
from the absorbing role of $0$, the neutrality of $e$, and the fact that
the product of two elements of $\{a,b\}$ is $0$.  Direct inspection against
the listed comparisons shows joint increasingness.  Hence $U_1$ is a
uninorm.  Nevertheless, $a\in\Ie$ is not right-transitive, since
$a\preceqT b\preceqT e$ but $a\npreceqT e$; therefore \textnormal{(C1)}
fails.  The same example shows that the lattice-style equality alone does
not imply all of \textnormal{(C1)--(C4)} on a proper trellis.
\end{example}

\begin{theorem}\label{thm:U2}\label{cor:U2-relative-characterization}
Let $\mathbb T$ be a bounded trellis, let $e\in X\setminus\{0,1\}$, and
let $\intr$ be an interior operator.  Let $T$ be a triangular norm on the
sub-trellis $[0,e]$ with neutral element $e$.  Assume \textnormal{(I1)} and
define
\begin{equation}\label{eq:U2}
U_2(x,y)=
\begin{cases}
T(x,y), & x,y\in[0,e],\\
x, & (x,y)\in\bigl((e,1]\times[0,e]\bigr)
      \cup\bigl(\Ie\times\{e\}\bigr),\\
y, & (x,y)\in\bigl([0,e]\times(e,1]\bigr)
      \cup\bigl(\{e\}\times\Ie\bigr),\\
\intr(x)\vee(x\wedge e), & x\in\Ie,\ y\in[0,e),\\
\intr(y)\vee(y\wedge e), & y\in\Ie,\ x\in[0,e),\\
(x\vee y)\vee e, & x,y\in\Ie,\\
1, & \text{otherwise}.
\end{cases}
\end{equation}
Then $U_2$ is a uninorm with neutral element $e$ if and only if
\textnormal{(I2)} and \textnormal{(I4)} hold.  Under these conditions,
\textnormal{(I3)} follows from \textnormal{(I1)} and \textnormal{(I2)}, and
$1$ is an absorbing element of $U_2$.
\end{theorem}

\begin{Proof}
Consider the order-dual bounded trellis $\mathbb T^{\partial}$.
With respect to $\mathbb T^{\partial}$, the interior operator $\intr$ is a
closure operator, the interval $[0,e]$ in $\mathbb T$ becomes $[e,1]$ in
$\mathbb T^{\partial}$, and $T$ becomes a triangular conorm with neutral
element $e$.  Conditions \textnormal{(I1)--(I4)} become, respectively,
\textnormal{(C1)--(C4)}, while $d_e$ and $n_e$ correspond to $c_e$ and $m_e$.
Formula \eqref{eq:U2} becomes
\eqref{eq:U1}.  The assertion therefore follows from \Cref{thm:U1} by order
duality.
\end{Proof}

\subsection{Variants with an absorbing endpoint}

\begin{theorem}\label{thm:U3}
Let $\mathbb T$ be a bounded trellis, let $e\in X\setminus\{0,1\}$, let
$\cl$ be a closure operator, and let $S$ be a triangular conorm on $[e,1]$
with neutral element $e$.  Assume \textnormal{(C1)--(C4)} and
$(e,1)\neq\varnothing$, and define $U_3\colon X^2\to X$ by
\begin{equation}\label{eq:U3}
U_{3}(x, y)=
\begin{cases}
1, & x=1 \text{ or } y=1, \\
S(x, y), & x, y \in[e, 1), \\
x, & (x,y)\in\bigl([0,e)\times[e,1)\bigr)
      \cup\bigl(\Ie\times\{e\}\bigr), \\
y, & (x,y)\in\bigl([e,1)\times[0,e)\bigr)
      \cup\bigl(\{e\}\times\Ie\bigr), \\
\cl(x)\wedge(x\vee e), & x\in\Ie,\ y\in(e,1), \\
\cl(y)\wedge(y\vee e), & y\in\Ie,\ x\in(e,1), \\
(x\wedge y)\wedge e, & x,y\in\Ie, \\
0, & \text{otherwise}.
\end{cases}
\end{equation}
Then $U_3$ is a uninorm on $\mathbb T$ with neutral element $e$ if and only if
\begin{equation}
S(x, y) \neq 1 \quad \text{for every } x, y \in(e, 1).
\label{eq:S_not_1}
\end{equation}
\end{theorem}

\begin{Proof}
We first prove necessity.  Suppose $S(x,y)=1$ for some $x,y\in(e,1)$.
Then
\[
 U_3\bigl(0,U_3(x,y)\bigr)=U_3(0,1)=1,
 \qquad
 U_3\bigl(U_3(0,x),y\bigr)=U_3(0,y)=0,
\]
so associativity fails.  Hence \eqref{eq:S_not_1} is necessary.

Conversely, assume \eqref{eq:S_not_1}.  Formula \eqref{eq:U3} is symmetric and
unambiguous, and direct substitution gives $U_3(e,x)=x$ for every $x\in X$.

To prove joint increasingness, let $x_1\preceqT x_2$ and
$y_1\preceqT y_2$.  If $x_2=1$ or $y_2=1$, then
$U_3(x_2,y_2)=1$, so the desired comparison follows from the maximality of
$1$.  Otherwise $x_2,y_2\ne1$.  Antisymmetry then implies
$x_1,y_1\ne1$: for example, $x_1=1$ together with
$1=x_1\preceqT x_2\preceqT1$ would force $x_2=1$.  On
$(X\setminus\{1\})^2$, the operations $U_3$ and $U_1$ agree pointwise.
Consequently, the joint increasingness proved in \Cref{thm:U1} gives
\[
 U_3(x_1,y_1)=U_1(x_1,y_1)
 \preceqT U_1(x_2,y_2)=U_3(x_2,y_2).
\]

We next show that $X\setminus\{1\}$ is closed under $U_3$.  For inputs in
$[e,1)$, condition \eqref{eq:S_not_1} applies when both inputs belong to
$(e,1)$; when one input equals $e$, neutral-element identity
$S(e,t)=t$ gives a value different from $1$.  A mixed incomparable--upper
clause has value $c_e(t)\in\Ie$ by (C2), an incomparable--incomparable
clause has value in $[0,e)$ by \Cref{lem:images}, the mixed clauses
return one of the inputs, and the remaining clause returns $0$. Hence
$X\setminus\{1\}$ is closed under $U_3$.

If at least one of $x,y,z$ equals $1$, both bracketings of their product are
$1$, because $1$ is absorbing for $U_3$.  If none equals $1$, the closure
just proved ensures that both intermediate products also differ from $1$.
Since $U_3=U_1$ on $(X\setminus\{1\})^2$, associativity of $U_1$ yields
\[
 \begin{aligned}
 U_3\bigl(U_3(x,y),z\bigr)
 &=U_1\bigl(U_1(x,y),z\bigr)\\
 &=U_1\bigl(x,U_1(y,z)\bigr)\\
 &=U_3\bigl(x,U_3(y,z)\bigr).
 \end{aligned}
\]
Therefore $U_3$ is a uninorm with neutral element $e$.
\end{Proof}

\begin{proposition}
\label{prop:U3-relative-characterization}
Let $\mathbb T$ be a bounded trellis, let
$e\in X\setminus\{0,1\}$, let $\cl$ be a closure operator, and let
$S$ be a triangular conorm on $[e,1]$ with neutral element $e$.
Assume \textnormal{(C1)}, and define $U_3$ by \eqref{eq:U3}.
\begin{enumerate}
\item If $(e,1)\ne\varnothing$, then the following statements are
equivalent:
\begin{enumerate}
\item $U_3$ is a uninorm on $\mathbb T$ with neutral element $e$;
\item conditions \textnormal{(C2)} and \textnormal{(C4)} hold, and
\begin{equation}\label{eq:U3-relative-endpoint}
 S(x,y)\ne1\qquad(x,y\in(e,1)).
\end{equation}
\end{enumerate}
Under these equivalent conditions, \textnormal{(C3)} follows from
\textnormal{(C1)} and \textnormal{(C2)} and is therefore redundant.

\item If $(e,1)=\varnothing$, then $U_3$ is a uninorm if and only if
$m_e$ is jointly increasing on $\Ie^2$.  In this degenerate case no
condition involving $c_e$ is necessary.
\end{enumerate}
\end{proposition}

\begin{Proof}
Assume first that $(e,1)\ne\varnothing$ and that $U_3$ is a uninorm.
Condition \eqref{eq:U3-relative-endpoint} is necessary.
Indeed, if
$p,q\in(e,1)$ and $S(p,q)=1$, then
\[
 U_3\bigl(0,U_3(p,q)\bigr)=U_3(0,1)=1,
 \qquad
 U_3\bigl(U_3(0,p),q\bigr)=U_3(0,q)=0,
\]
contrary to associativity.  Hence
\eqref{eq:U3-relative-endpoint} holds.

Condition \textnormal{(C1)} also implies
\begin{equation}\label{eq:U3-me-lower}
 m_e(x,y)\in[0,e)\qquad(x,y\in\Ie).
\end{equation}
Indeed, $m_e(x,y)\preceqT e$.  If $m_e(x,y)=e$, then
$e\preceqT x\wedge y\preceqT x$,
and right-transitivity of $e$ would give $e\preceqT x$, contrary to
$x\parallel e$.

Fix $a\in(e,1)$ and $x\in\Ie$.  Since $x\preceqT c_e(x)$, if
$c_e(x)\preceqT e$, right-transitivity of $x$ would imply
$x\preceqT e$, a contradiction.  Thus $c_e(x)\in\Ie\cup(e,1]$.
Associativity at the triple $(x,x,a)$ gives
\begin{equation}\label{eq:U3-image-test}
 \begin{aligned}
 m_e(x,x)
 &=U_3\bigl(U_3(x,x),a\bigr)\\
 &=U_3\bigl(x,U_3(x,a)\bigr)
 =U_3\bigl(x,c_e(x)\bigr).
 \end{aligned}
\end{equation}
If $c_e(x)\in(e,1]$, then the last member of
\eqref{eq:U3-image-test} equals $c_e(x)$, whereas its first member
belongs to $[0,e)$ by \eqref{eq:U3-me-lower}.  This is impossible.
Consequently,
\begin{equation}\label{eq:U3-ce-image}
 c_e(\Ie)\subseteq\Ie.
\end{equation}

If $x,y\in\Ie$ and $x\preceqT y$, increasingness of $U_3$, applied
to $(x,a)\preceqT(y,a)$, gives
$c_e(x)=U_3(x,a)\preceqT U_3(y,a)=c_e(y)$.
Hence $c_e$ is increasing on $\Ie$.

Increasingness of $S$ and neutrality of $e$ yield
$a=S(e,a)\preceqT S(a,a)$.
Since $a\ne e$, antisymmetry shows that $S(a,a)\ne e$.
Together with \eqref{eq:U3-relative-endpoint}, this gives
$S(a,a)\in(e,1)$.  Associativity at $(x,a,a)$ and
\eqref{eq:U3-ce-image} now yield
\[
 \begin{aligned}
 c_e(c_e(x))
 &=U_3\bigl(U_3(x,a),a\bigr)\\
 &=U_3\bigl(x,U_3(a,a)\bigr)\\
 &=U_3\bigl(x,S(a,a)\bigr)\\
 &=c_e(x).
 \end{aligned}
\]
Thus \textnormal{(C2)} holds.

The restriction of $U_3$ to $\Ie^2$ is $m_e$.  Therefore the joint
increasingness of $U_3$ implies that $m_e$ is jointly increasing on
$\Ie^2$.  Finally, for $x,y\in\Ie$, associativity at $(a,x,y)$ gives
\[
 \begin{aligned}
 m_e(c_e(x),y)
 &=U_3\bigl(U_3(a,x),y\bigr)\\
 &=U_3\bigl(a,U_3(x,y)\bigr)\\
 &=U_3\bigl(a,m_e(x,y)\bigr)\\
 &=m_e(x,y),
 \end{aligned}
\]
where \eqref{eq:U3-me-lower} is used in the last equality.  Hence
\textnormal{(C4)} holds.

Conversely, suppose that \textnormal{(C2)}, \textnormal{(C4)}, and
\eqref{eq:U3-relative-endpoint} hold.  As in
\Cref{prop:U1-relative-characterization}, conditions \textnormal{(C1)}
and \textnormal{(C2)} imply \textnormal{(C3)}.  Therefore
\textnormal{(C1)--(C4)} hold, and \Cref{thm:U3} shows that $U_3$ is a
uninorm.

It remains to consider $(e,1)=\varnothing$.  In this case the
definition of $U_3$ reduces to
\[
 U_3(x,y)=
 \begin{cases}
 1, & x=1\ \text{or}\ y=1,\\
 y, & x=e,\\
 x, & y=e,\\
 m_e(x,y), & x,y\in\Ie,\\
 0, & \text{otherwise}.
 \end{cases}
\]
In particular, neither $c_e$ nor a nontrivial value of $S$ occurs in
the formula.

If $U_3$ is a uninorm, its restriction to $\Ie^2$ shows immediately
that $m_e$ is jointly increasing.  Conversely, assume that $m_e$ is
jointly increasing.  Commutativity and neutrality of $e$ follow
directly from the displayed formula.

For increasingness, let $x\preceqT x'$ and $y\preceqT y'$.  If
$x'=1$ or $y'=1$, the required comparison follows from the maximality
of $1$.  Otherwise, \Cref{lem:transitions} and
$(e,1)=\varnothing$ imply that $x=e$ forces $x'=e$, and
$x\in\Ie$ forces $x'\in\Ie$; the same statements hold for $y$ and
$y'$.  If $x=e$ or $y=e$, neutrality reduces the desired inequality to
the comparison involving the other variable.  If $x,y\in\Ie$, the
assertion is the joint increasingness of $m_e$.  In every remaining case
$U_3(x,y)=0$.  Thus $U_3$ is jointly increasing.

For associativity, a triple containing $1$ gives the value $1$ under
both bracketings, and a triple containing $e$ reduces by neutrality.
Every remaining entry belongs to $[0,e)\cup\Ie$.  If all three
entries belong to $\Ie$, each inner product belongs to $[0,e)$ by
\eqref{eq:U3-me-lower}, and the second product is $0$.  If exactly two
entries belong to $\Ie$, pairing them first gives an element of
$[0,e)$ and hence gives $0$ at the second step, while pairing an
incomparable and a lower element first gives $0$ immediately.  If at
most one entry belongs to $\Ie$, both bracketings are again $0$.
Therefore $U_3$ is associative, completing the proof.
\end{Proof}

\begin{theorem}\label{thm:U4}
Let $\mathbb T$ be a bounded trellis, let $e\in X\setminus\{0,1\}$, let
$\intr$ be an interior operator, and let $T$ be a triangular norm on $[0,e]$
with neutral element $e$.  Assume \textnormal{(I1)--(I4)} and
$(0,e)\ne\varnothing$, and define
\begin{equation}\label{eq:U4}
U_4(x,y)=
\begin{cases}
0, & x=0\ \text{or}\ y=0,\\
T(x,y), & x,y\in(0,e],\\
x, & (x,y)\in\bigl((e,1]\times(0,e]\bigr)
      \cup\bigl(\Ie\times\{e\}\bigr),\\
y, & (x,y)\in\bigl((0,e]\times(e,1]\bigr)
      \cup\bigl(\{e\}\times\Ie\bigr),\\
\intr(x)\vee(x\wedge e), & x\in\Ie,\ y\in(0,e),\\
\intr(y)\vee(y\wedge e), & y\in\Ie,\ x\in(0,e),\\
(x\vee y)\vee e, & x,y\in\Ie,\\
1, & \text{otherwise}.
\end{cases}
\end{equation}
Then $U_4$ is a uninorm if and only if
\begin{equation}\label{eq:T_not_0}
 T(x,y)\ne0\qquad\text{for every }x,y\in(0,e).
\end{equation}
\end{theorem}

\begin{Proof}
We first prove necessity.  If $T(x,y)=0$ for some $x,y\in(0,e)$, then
\[
 U_4\bigl(1,U_4(x,y)\bigr)=U_4(1,0)=0,
 \qquad
 U_4\bigl(U_4(1,x),y\bigr)=U_4(1,y)=1.
\]
Thus associativity implies \eqref{eq:T_not_0}.

Conversely, assume \eqref{eq:T_not_0}.  The formula is the order dual of the definition
in \Cref{thm:U3}, so it is well defined and symmetric, and $e$ is neutral.
We give the remaining dual argument explicitly.

Let $x_1\preceqT x_2$ and $y_1\preceqT y_2$.  If $x_1=0$ or $y_1=0$,
then $U_4(x_1,y_1)=0\preceqT U_4(x_2,y_2)$.  Otherwise
$x_1,y_1\ne0$, and antisymmetry implies $x_2,y_2\ne0$; for example,
$0\preceqT x_1\preceqT x_2=0$ would force $x_1=0$.  On
$(X\setminus\{0\})^2$, the operations $U_4$ and $U_2$ agree pointwise.
The joint increasingness of $U_2$ therefore gives
\[
 U_4(x_1,y_1)=U_2(x_1,y_1)
 \preceqT U_2(x_2,y_2)=U_4(x_2,y_2).
\]

The set $X\setminus\{0\}$ is closed under $U_4$.  Indeed, for inputs in
$(0,e]$, condition \eqref{eq:T_not_0} treats the case in which both inputs
belong to $(0,e)$, while $T(e,t)=t$ treats the case in which one input is
$e$.  The mixed incomparable--lower clauses take values in $\Ie$ by (I2),
the incomparable--incomparable clause takes values in $(e,1]$ by
\Cref{lem:images}, the mixed clauses return a nonzero input, and the
remaining clause returns $1$.

If one of $x,y,z$ is $0$, both bracketings are $0$, because $0$ is absorbing
for $U_4$.  If none is $0$, closure keeps both intermediate values in
$X\setminus\{0\}$, where $U_4=U_2$.  Hence
\[
 \begin{aligned}
 U_4\bigl(U_4(x,y),z\bigr)
 &=U_2\bigl(U_2(x,y),z\bigr)\\
 &=U_2\bigl(x,U_2(y,z)\bigr)\\
 &=U_4\bigl(x,U_4(y,z)\bigr)
 \end{aligned}
\]
by associativity of $U_2$.  Thus $U_4$ is a uninorm.
\end{Proof}

Combining the preceding result with the order-dual necessity argument gives
the criterion for the interior construction with absorbing element $0$.

\begin{corollary}
\label{cor:U4-relative-characterization}
Let $\mathbb T$ be a bounded trellis,
$e\in X\setminus\{0,1\}$,   $\intr$ be an interior operator, and
$T$ be a triangular norm on $[0,e]$ with neutral element $e$.
Assume \textnormal{(I1)}, and define $U_4$ by \eqref{eq:U4}.
\begin{enumerate}
\item If $(0,e)\ne\varnothing$, then $U_4$ is a uninorm if and only if
\textnormal{(I2)}, \textnormal{(I4)}, and $T(x,y)\ne0$ for all
$x,y\in(0,e)$ hold.  Under these conditions, \textnormal{(I3)} is redundant.
\item If $(0,e)=\varnothing$, then $U_4$ is a uninorm if and only if
$n_e$ is jointly increasing on $\Ie^2$; no condition involving $d_e$
is necessary.
\end{enumerate}
\end{corollary}

\begin{Proof}
Apply \Cref{prop:U3-relative-characterization} to the order-dual bounded
trellis. Under this duality, (C1)--(C4), $c_e$, $m_e$, and
condition \eqref{eq:U3-relative-endpoint} become
(I1)--(I4), $d_e$, $n_e$, and
condition \eqref{eq:T_not_0}, respectively. The nondegenerate and degenerate
criteria therefore dualize exactly as stated.
\end{Proof}

\subsection{Constructions using both a triangular norm and a triangular conorm}

Let $\mathbb T$ be a bounded trellis, let $e\in X\setminus\{0,1\}$, let
$\cl$ be a closure operator, let $T$ be a triangular norm on $[0,e]$ with
neutral element $e$, and let $S$ be a triangular conorm on $[e,1]$ with
neutral element $e$.  We shall use the following sufficient condition:
\begin{equation}\label{eq:separation-lower}
 x\parallel y\qquad\text{for every }x\in(0,e),\ y\in\Ie.
\end{equation}
Define
\begin{equation}\label{eq:U5}
U_5(x,y)=
\begin{cases}
T(x,y), & x,y\in[0,e),\\
S(x,y), & x,y\in[e,1],\\
x, & (x,y)\in\bigl([0,e)\times[e,1]\bigr)
      \cup\bigl(\Ie\times\{e\}\bigr),\\
y, & (x,y)\in\bigl([e,1]\times[0,e)\bigr)
      \cup\bigl(\{e\}\times\Ie\bigr),\\
\cl(x)\wedge(x\vee e), & x\in\Ie,\ y\in(e,1],\\
\cl(y)\wedge(y\vee e), & y\in\Ie,\ x\in(e,1],\\
0, & \text{otherwise}.
\end{cases}
\end{equation}

The following result gives the necessary and sufficient condition under
\textnormal{(C1)}.

\begin{theorem}\label{thm:U5}\label{prop:U5-relative-characterization}
Let $\mathbb T$ be a bounded trellis, let $e\in X\setminus\{0,1\}$, let
$\cl$ be a closure operator, let $T$ be a triangular norm on $[0,e]$ with
neutral element $e$, and let $S$ be a triangular conorm on $[e,1]$ with
neutral element $e$.  Assume (C1), and define $U_5$ by \eqref{eq:U5}.
Then the following statements are equivalent:
\begin{enumerate}
\item $U_5$ is a uninorm on $\mathbb T$ with neutral element $e$;
\item condition (C2) holds and
\begin{equation}\label{eq:zero-comparison-lower}
 p\preceqT z\quad\Longrightarrow\quad T(p,q)=0
 \qquad\bigl(p,q\in[0,e),\ z\in\Ie\bigr).
\end{equation}
\end{enumerate}
Under these assumptions, (C3) follows from (C2), whereas (C4) is not
required for this construction.
\end{theorem}

\begin{Proof}
Suppose first that $U_5$ is a uninorm.  For $x\in\Ie$, put
$c=c_e(x)$.  By \Cref{lem:images}, $x\preceqT c$.  If
$c\preceqT e$, then right-transitivity of $x$ would imply
$x\preceqT e$, contrary to $x\parallel e$.  Hence
$c\in\Ie\cup(e,1]$.

Assume that $c\in(e,1]$.  Since $1$ is absorbing for the triangular
conorm $S$, associativity at $(x,1,1)$ gives
\[
 1=U_5(c,1)
  =U_5\bigl(U_5(x,1),1\bigr)
  =U_5\bigl(x,U_5(1,1)\bigr)
  =c.
\]
Thus $c=1$.  Associativity at $(x,x,1)$ would then give
\[
 0=U_5\bigl(U_5(x,x),1\bigr)
  =U_5\bigl(x,U_5(x,1)\bigr)=1,
\]
a contradiction.  Consequently $c_e(\Ie)\subseteq\Ie$.
If $x\preceqT y$ in $\Ie$, joint increasingness applied to
$(x,1)\preceqT(y,1)$ gives
$c_e(x)=U_5(x,1)\preceqT U_5(y,1)=c_e(y)$.
Moreover, associativity at $(x,1,1)$, now using
$c_e(x)\in\Ie$, yields
\[
 c_e(c_e(x))
 =U_5\bigl(U_5(x,1),1\bigr)
 =U_5\bigl(x,U_5(1,1)\bigr)=c_e(x).
\]
This proves (C2).

Let $p,q\in[0,e)$ and $z\in\Ie$ satisfy $p\preceqT z$.  Then
$(p,q)\preceqT(z,q)$ componentwise, by the definition of $U_5$,
$U_5(z,q)=0$ for $z\in\Ie$ and $q\in [0,e)$, and joint increasingness gives
$T(p,q)=U_5(p,q)\preceqT U_5(z,q)=0$.
Since $0$ is the least element, \eqref{eq:zero-comparison-lower} follows.

Conversely, assume (C2) and \eqref{eq:zero-comparison-lower}.  We first
recover (C3).  If $x\in\Ie$, $y\in(e,1]$, and $x\preceqT y$, then $y$
is an upper bound of both $x$ and $e$, and hence
$c_e(x)\preceqT x\vee e\preceqT y$.
Condition (C2) gives $c_e(x)\in\Ie$, so (C1) makes $c_e(x)$
right-transitive.  Therefore $c_e(x)\preceqT y$, which is (C3).

We next prove joint increasingness.  If $p,q\in[0,e)$, monotonicity and
neutrality of $T$ give
\[
 t=T(p,q)\preceqT T(p,e)=p,
 \qquad
 t\preceqT T(e,q)=q,
\]
and $t\in[0,e)$.  Similarly, for $r,s\in(e,1]$ one has
$r,s\preceqT S(r,s)$ and $S(r,s)\in(e,1]$.
Let $x\preceqT x'$ and $y\preceqT y'$.  By commutativity, we may
simultaneously interchange $x$ with $y$ and $x'$ with $y'$ in the cases
below.

Suppose first that $x,y\in[0,e)$ and put $t=T(x,y)$.  If
$x',y'\in[0,e)$, the assertion follows from joint increasingness of $T$.
If one of $x',y'$ lies in $[0,e)$ and the other in $\Ie$, or if both lie in
$\Ie$, then $x\preceqT x'$ or $y\preceqT y'$ has an element of $\Ie$ on
the right.  Condition \eqref{eq:zero-comparison-lower}, together with
commutativity of $T$, gives $t=0=U_5(x',y')$.  If, say,
$x'\in[0,e)$ and $y'\in[e,1]$, then $U_5(x',y')=x'$ and
$t\preceqT x\preceqT x'$; right-transitivity of $t$ completes the
comparison.  If $x',y'\in[e,1]$, then
$t\preceqT e\preceqT S(x',y')$,
and middle-transitivity of $e$ applies.  If $\{x',y'\}=\{e,z\}$ with
$z\in\Ie$, right-transitivity of $t$ gives
$t\preceqT z=U_5(x',y')$.  If one of $x',y'$ is $z\in\Ie$ and the other
lies in $(e,1]$, the same argument gives
$t\preceqT z\preceqT c_e(z)=U_5(x',y')$.

Suppose next that $x=p\in[0,e)$ and $y\in[e,1]$.  Then
$U_5(x,y)=p$ and \Cref{lem:transitions} gives $y'\in[e,1]$.  If
$x'\in[0,e)$, then $U_5(x',y')=x'$ and $p\preceqT x'$.  If
$x'\in[e,1]$, then $p\preceqT e\preceqT U_5(x',y')$, and
middle-transitivity of $e$ applies.  If $x'=z\in\Ie$, then
$U_5(x',y')$ is $z$ when $y'=e$ and is $c_e(z)$ when $y'\in(e,1]$; the
required comparison follows from $p\preceqT z\preceqT c_e(z)$ and
right-transitivity of $p$.  If $x\in[0,e)$ and $y\in\Ie$, then
$U_5(x,y)=0$.

If $x,y\in[e,1]$, then $x',y'\in[e,1]$ and joint increasingness of $S$
applies.  Suppose next that $x=e$ and $y=z\in\Ie$.  If
$y'=z'\in\Ie$, then $U_5(x',y')$ is either $z'$ or $c_e(z')$, and
$z\preceqT z'\preceqT c_e(z')$ together with right-transitivity of $z$
gives the result.  If $y'\in(e,1]$, then $U_5(x',y')$ is either $y'$ or
an $S$-value above $y'$, and the same argument applies.

Finally, suppose that $x\in(e,1]$ and $y=z\in\Ie$.  If $y'\in\Ie$,
increasingness of $c_e$ applies.  If $y'\in(e,1]$, condition (C3) gives
$c_e(z)\preceqT y'\preceqT S(x',y')$,
after a possible simultaneous interchange of the coordinates; since
$c_e(z)\in\Ie$ is right-transitive, the comparison follows.  If
$x,y\in\Ie$, then $U_5(x,y)=0$.  These cases prove joint increasingness of
$U_5$ without using \eqref{eq:separation-lower} or (C4).

It remains to prove associativity.  Triples containing $e$ reduce by
neutrality.  On $[0,e)^3$ and $(e,1]^3$, associativity follows from that of
$T$ and $S$, respectively, while every pairing on $\Ie^3$ is $0$.  For two
lower and one upper input, every pairing gives their $T$-value; for two
lower and one incomparable input, every pairing gives $0$.  For one lower
and two upper inputs, every pairing gives the lower input, and for one lower
and two incomparable inputs every pairing gives $0$, using $T(0,x)=0$.
For one upper and two incomparable inputs, (C2) places the two possible
$c_e$-values in $\Ie$, so all pairings are $0$.  For one incomparable and
two upper inputs, all pairings give $c_e(x)$ by the idempotence in (C2).
For one input from each of $[0,e)$, $(e,1]$, and $\Ie$, every pairing is
$0$.  Thus $U_5$ is associative.  Commutativity and neutrality of $e$ are
immediate from \eqref{eq:U5}, completing the proof.
\end{Proof}

\begin{corollary}\label{cor:U5-separation}
Assume \textnormal{(C1)--(C3)} and \eqref{eq:separation-lower}.  Then the
operation $U_5$ in \eqref{eq:U5} is a uninorm with neutral element $e$.
\end{corollary}

\begin{Proof}
Let $p,q\in[0,e)$ and $z\in\Ie$ with $p\preceqT z$.  Condition
\eqref{eq:separation-lower} excludes $p\in(0,e)$, so $p=0$ and
$T(p,q)=0$.  Hence \eqref{eq:zero-comparison-lower} holds, and the result
follows from \Cref{thm:U5}.
\end{Proof}

The incomparability restrictions used in the closure-based trellis
constructions of \cite{KongZhao2024} are sufficient for controlling pairs
consisting of a lower and an incomparable element.  For the fixed formula
\eqref{eq:U5}, \Cref{cor:U5-separation} gives the corresponding sufficient
statement, whereas \Cref{thm:U5} shows that the exact requirement is
\eqref{eq:zero-comparison-lower}.  The two requirements are equivalent when
$T$ has no zero divisors by \Cref{cor:U5-no-zero-divisors};
\Cref{rem:cond9_sufficient_not_necessary} shows that they can differ when zero
divisors are present.  Thus the improvement concerns the criterion for
\eqref{eq:U5}, rather than an inclusion between all classes constructed in the
two papers.

The construction in \Cref{cor:U5-separation} occurs on the following small
proper bounded trellis.

\begin{example}\label{ex:proper}
Let $X=\{0,e,d,a,b,c,1\}$.  Besides reflexivity and
$0\preceqT x\preceqT1$, impose
$e\precT a$, $e\precT b$, $e\precT c$ and
$a\precT b\precT c\precT a$, with no further comparisons.  Comparable
pairs have the indicated lower and upper elements, whereas every
incomparable pair has meet $0$ and join $1$.  Hence $X$ is a bounded
trellis.  It is proper because $a\preceqT b\preceqT c$ but
$a\npreceqT c$.

Define $\cl(a)=\cl(b)=\cl(c)=1$ and let $\cl$ fix $0,e,d,1$.
This map is extensive and idempotent.  It is also increasing: the only
nonboundary comparisons not fixed by $\cl$ have an element of
$\{a,b,c\}$ on the right, and all three of these elements have image $1$.
Thus $\cl$ is a closure operator.  Moreover,
$\Ie=\{d\}$, $[0,e)=\{0\}$, and $(e,1]=\{a,b,c,1\}$.
The element $e$ is middle- and right-transitive, while $0$ and $d$ are
right-transitive, so \textnormal{(C1)} holds.  Moreover,
$c_e(d)=d$ and $m_e(d,d)=0$; since $\Ie$ is a singleton, conditions
\textnormal{(C2)--(C4)} follow immediately.  Condition
\eqref{eq:separation-lower} is vacuous because $(0,e)=\varnothing$; hence
\eqref{eq:zero-comparison-lower} holds.

On $[0,e]=\{0,e\}$ take the unique triangular norm.  On
$[e,1]=\{e,a,b,c,1\}$ define $S(e,x)=S(x,e)=x$ and put $S(x,y)=1$
whenever $x,y\ne e$.  The operation is commutative and has neutral element
$e$.  If at least two of three inputs differ from $e$, both associative
products are $1$; otherwise associativity follows from neutrality.  For
increasingness, a value different from $1$ can occur only when one input is
$e$, in which case the required comparison is one of the assumed coordinate
comparisons.  Hence $S$ is a triangular conorm, and \Cref{thm:U5} yields a
uninorm on this proper bounded trellis.
\end{example}

\begin{corollary}\label{cor:U5-no-zero-divisors}
Under the hypotheses of \Cref{prop:U5-relative-characterization}, suppose
that $T$ has no zero divisors, that is,
\[
 T(p,q)=0\quad\Longrightarrow\quad p=0\ \text{or}\ q=0
 \qquad(p,q\in[0,e]).
\]
Then \eqref{eq:zero-comparison-lower} is equivalent to
\eqref{eq:separation-lower}.  Consequently, under (C1) and the
no-zero-divisor assumption, condition (C2) together with
\eqref{eq:separation-lower} is necessary and sufficient for $U_5$ to be a
uninorm.
\end{corollary}

\begin{Proof}
Assume first \eqref{eq:separation-lower}.  If $p\in[0,e)$ and
$p\preceqT z\in\Ie$, then $p$ cannot belong to $(0,e)$, so $p=0$; hence
$T(p,q)=T(0,q)=0$.  Thus \eqref{eq:zero-comparison-lower} holds.
Conversely, assume \eqref{eq:zero-comparison-lower}.  If
$p\in(0,e)$ and $p\preceqT z\in\Ie$, taking $q=p$ gives $T(p,p)=0$, which
contradicts the absence of zero divisors.  The reverse comparison
$z\preceqT p$ is also impossible: since $z\in\Ie\subseteq X^{\rtr}$,
the chain $z\preceqT p\preceqT e$ would imply $z\preceqT e$.  Hence
$p\parallel z$, as required.
\end{Proof}

The converse relation between \eqref{eq:zero-comparison-lower} and
\eqref{eq:separation-lower} is clarified by the following counterexample.

\begin{remark}
\label{rem:cond9_sufficient_not_necessary}
Condition \eqref{eq:separation-lower} always implies
\eqref{eq:zero-comparison-lower}, but the converse may fail when $T$ has
zero divisors.  Indeed, let
$
 X=\{0,p,e,z,1\}
$
be the bounded lattice whose nontrivial covering relations are
\[
 0\precT p\precT e\precT1,
 \qquad p\precT z\precT1,
 \qquad e\parallel z.
\]
Take $\cl$ to be the identity map.  Then $\Ie=\{z\}$ and
$c_e(z)=z\wedge(z\vee e)=z\wedge1=z$,
so (C1) and (C2) hold.  On the three-element chain
$[0,e]=\{0,p,e\}$ define the drastic triangular norm
\[
 T(u,v)=
 \begin{cases}
 u\wedge v,&u=e\ \text{or}\ v=e,\\
 0,&u,v\in\{0,p\},
 \end{cases}
\]
and take the triangular conorm on $[e,1]=\{e,1\}$ determined by
$S(e,e)=e$ and $S(u,v)=1$ otherwise.  Since
$[0,e)=\{0,p\}$, one has $T(u,v)=0$ for all $u,v\in[0,e)$,
so \eqref{eq:zero-comparison-lower} holds.  Hence
\Cref{prop:U5-relative-characterization} shows that $U_5$ is a uninorm.
Nevertheless, $p\preceqT z$, and therefore
\eqref{eq:separation-lower} fails.  Thus condition
\eqref{eq:separation-lower} is sufficient but is not necessary in general.
\end{remark}

\begin{theorem}\label{thm:U6}
Let $\intr$ be an interior operator,   $T$ be a triangular norm on
$[0,e]$, and   $S$ be a triangular conorm on $[e,1]$, both with neutral
element $e$.  Assume \textnormal{(I1)--(I3)} and the condition
\begin{equation}\label{eq:separation-upper}
 x\parallel y\qquad\text{for every }x\in(e,1),\ y\in\Ie.
\end{equation}
Define
\begin{equation}\label{eq:U6}
U_6(x,y)=
\begin{cases}
T(x,y), & x,y\in[0,e],\\
S(x,y), & x,y\in(e,1],\\
x, & (x,y)\in\bigl((e,1]\times[0,e]\bigr)
      \cup\bigl(\Ie\times\{e\}\bigr),\\
y, & (x,y)\in\bigl([0,e]\times(e,1]\bigr)
      \cup\bigl(\{e\}\times\Ie\bigr),\\
\intr(x)\vee(x\wedge e), & x\in\Ie,\ y\in[0,e),\\
\intr(y)\vee(y\wedge e), & y\in\Ie,\ x\in[0,e),\\
1, & \text{otherwise}.
\end{cases}
\end{equation}
Then $U_6$ is a uninorm with neutral element $e$.
\end{theorem}

\begin{Proof}
Consider the order-dual bounded trellis $\mathbb T^{\partial}$.  The
incomparable set of $e$ is unchanged.

With respect to $\mathbb T^{\partial}$, the map $\intr$ is a closure
operator.  Indeed, $\intr$ is extensive in the dual order because
$\intr(x)\preceqT x$ in the original order implies
$x\preceqT^{\partial}\intr(x)$.  Conditions (I1)--(I3) are exactly
conditions (C1)--(C3) in the dual pseudo-order, and the expression
corresponding to $c_e$ is $\intr(x)\vee(x\wedge e)=d_e(x)$.
The triangular conorm $S$ on the original interval $[e,1]$ is a triangular
norm on the dual interval from $1$ to $e$, whereas $T$ is the corresponding
triangular conorm on the dual interval from $e$ to $0$.  Moreover,
\eqref{eq:separation-upper} is precisely
\eqref{eq:separation-lower} in the dual pseudo-order.

Applying \Cref{cor:U5-separation} to these dual data gives exactly the cases displayed
in \eqref{eq:U6}: the lower and upper intervals are interchanged, $0$ and
$1$ are interchanged, and $c_e$ is replaced by $d_e$.  Hence the joint
increasingness and all three associativity pairings proved for $U_5$ dualize
without any additional assumption.  Therefore $U_6$ is a uninorm with
neutral element $e$.
\end{Proof}

We next state the exact criterion for the two-sided interior construction.

\begin{corollary}
\label{cor:U6-relative-characterization}
Let $\mathbb T$ be a bounded trellis,   $e\in X\setminus\{0,1\}$,
$\intr$ be an interior operator,   $T$ be a triangular norm on $[0,e]$
with neutral element $e$, and   $S$ be a triangular conorm on $[e,1]$
with neutral element $e$.  Assume (I1), and define $U_6$ by
\eqref{eq:U6}.  Then $U_6$ is a uninorm if and only if (I2) holds and
\begin{equation}\label{eq:one-comparison-upper}
 z\preceqT p\quad\Longrightarrow\quad S(p,q)=1
 \qquad\bigl(p,q\in(e,1],\ z\in\Ie\bigr).
\end{equation}
Under these assumptions, (I3) follows from (I2), whereas (I4) is not
required.  If $S$ has no nontrivial one-divisors, in the sense that
\[
 S(p,q)=1\quad\Longrightarrow\quad p=1\ \text{or}\ q=1
 \qquad(p,q\in[e,1]),
\]
then \eqref{eq:one-comparison-upper} is equivalent to
\eqref{eq:separation-upper}.  Consequently, under (I1) and the
absence of nontrivial one-divisors, condition (I2) together with
\eqref{eq:separation-upper} is necessary and sufficient for $U_6$ to be a
uninorm.
\end{corollary}

\begin{Proof}
Apply \Cref{prop:U5-relative-characterization} and
\Cref{cor:U5-no-zero-divisors} to the order-dual bounded trellis. Under this
duality, (C1)--(C3), condition
\eqref{eq:zero-comparison-lower}, and the absence of zero divisors become
(I1)--(I3), condition
\eqref{eq:one-comparison-upper}, and the absence of
nontrivial one-divisors, respectively. This gives all the stated equivalences.
\end{Proof}

\section{BOUNDED-LATTICE REDUCTIONS AND SCOPE}\label{sec:reduction}

Conditions \textnormal{(C1)--(C4)} and \textnormal{(I1)--(I4)} are
formulated for trellises because ordinary lattice manipulations cannot be
used without transitivity and associativity.  On a bounded lattice, these
conditions reduce to explicit identities and restrictions on the operator
images.  Combined with the characterizations in
\Cref{sec:constructions}, these reductions yield exact lattice criteria for
all six constructions.

\begin{proposition}\label{prop:lattice-reduction}
Let $L$ be a bounded lattice and let $e\in L\setminus\{0,1\}$.
\begin{enumerate}
\item If $\cl$ is a closure operator, then \textnormal{(C1)--(C4)} hold if
and only if $x\wedge e=\cl(x)\wedge e$ for all $x\in\Ie$.
\item If $\intr$ is an interior operator, then \textnormal{(I1)--(I4)} hold
if and only if $x\vee e=\intr(x)\vee e$ for all $x\in\Ie$.
\item If $\cl$ is a closure operator, then the following statements are
equivalent:
\begin{enumerate}
\item conditions \textnormal{(C1)--(C3)} hold;
\item $\cl(x)\parallel e$ for every $x\in\Ie$;
\item $\cl(x)\notin[e,1]$ for every $x\in\Ie$.
\end{enumerate}
Dually, for an interior operator $\intr$, one has
\[
\begin{aligned}
 \textnormal{(I1)--(I3) hold}
 &\quad\Longleftrightarrow\quad
 \intr(x)\parallel e &&(x\in\Ie)\\
 &\quad\Longleftrightarrow\quad
 \intr(x)\notin[0,e] &&(x\in\Ie).
\end{aligned}
\]
\end{enumerate}
\end{proposition}

\begin{Proof}
Every element of a lattice is right-, left-, and middle-transitive, and lattice meet and join are associative and jointly increasing.  Thus (C1) and its order dual hold automatically.

For (1), let $x\in\Ie$ and put $c_e(x)=\cl(x)\wedge(x\vee e)$.  Since $x\le c_e(x)$, the inequality $c_e(x)\le e$ would imply $x\le e$.  If $e\le c_e(x)$, then
\[
 e=c_e(x)\wedge e
  =\cl(x)\wedge(x\vee e)\wedge e
  =\cl(x)\wedge e=x\wedge e,
\]
which would imply $e\le x$.  Both alternatives contradict $x\parallel e$, so
$c_e(x)\in\Ie$.  Monotonicity of $c_e$ follows from that of $\cl$, meet, and
join.  Moreover, $x\le c_e(x)\le\cl(x)$ and
$x\le c_e(x)\le x\vee e$,
whence $\cl(c_e(x))=\cl(x)$ and $c_e(x)\vee e=x\vee e$.  Therefore $c_e(c_e(x))=c_e(x)$, proving (C2).  Joint increasingness of $m_e$ follows from that of meet, and
\[
\begin{aligned}
 m_e(c_e(x),y)
 &=\cl(x)\wedge(x\vee e)\wedge y\wedge e\\
 &=x\wedge y\wedge e=m_e(x,y),
\end{aligned}
\]
so (C4) holds.  Finally, if $x\le y$ and $y>e$, then $c_e(x)\le x\vee e\le y$, which is (C3).  This proves the sufficiency direction in (1).

Conversely, suppose that \textnormal{(C1)--(C4)} hold.  For $x\in\Ie$ put
$c=c_e(x)$.  Condition (C2) gives $c\in\Ie$, so (C4) may be applied with
$y=c$.  Since $x\le c$, lattice absorption and associativity of meet give
$c\wedge e=m_e(c,c)=m_e(x,c)=x\wedge e$.
On the other hand,
\[
 c\wedge e
 =\bigl[\cl(x)\wedge(x\vee e)\bigr]\wedge e
 =\cl(x)\wedge e,
\]
because $e\le x\vee e$.  Hence $x\wedge e=\cl(x)\wedge e$, proving the
necessity direction.  Assertion~(2), including both directions, is its
order dual.

For (3), fix $x\in\Ie$.  Extensivity gives $x\le\cl(x)$, so
$\cl(x)\le e$ would imply $x\le e$, which is impossible.  Consequently,
$\cl(x)$ can be comparable with $e$ only in the direction
$e\le\cl(x)$.  Since every element is below $1$, this proves the
equivalence $\cl(x)\parallel e\Longleftrightarrow\cl(x)\notin[e,1]$.

Suppose that $\cl(x)\parallel e$ for every $x\in\Ie$, and put
$c=c_e(x)$.  As before, $x\le c$, so $c\le e$ is impossible.  If
$e\le c$, then $c\le\cl(x)$ would give $e\le\cl(x)$, again a
contradiction.  Hence $c\in\Ie$.  The map $c_e$ is increasing because
closure, meet, and join are increasing on a lattice.  Moreover,
$x\le c_e(x)\le\cl(x)$ and $x\le c_e(x)\le x\vee e$, so closure
idempotence and lattice absorption give $\cl(c_e(x))=\cl(x)$ and
$c_e(x)\vee e=x\vee e$.
It follows that $c_e(c_e(x))=c_e(x)$, proving (C2).  If
$x\le y$ with $y>e$, then $y$ is an upper bound of both $x$ and $e$;
thus $c_e(x)\le x\vee e\le y$,
which is (C3).  Condition (C1) is automatic on a lattice, and hence
\textnormal{(C1)--(C3)} hold.

Conversely, suppose that \textnormal{(C1)--(C3)} hold.  By (C2),
$c_e(x)\in\Ie$.  If $e\le\cl(x)$, then $e$ is a lower bound of both
$\cl(x)$ and $x\vee e$, whence
$e\le\cl(x)\wedge(x\vee e)=c_e(x)$,
contrary to $c_e(x)\parallel e$.  Therefore $e\nleq\cl(x)$, and the
first paragraph shows that $\cl(x)\parallel e$, or equivalently
$\cl(x)\notin[e,1]$.  This completes the closure part of (3).  The
interior statement follows by order duality.
\end{Proof}

\begin{corollary}\label{cor:lattice_case_U1}
Let $L$ be a bounded lattice, let $e\in L\setminus\{0,1\}$, let $S$ be a
triangular conorm on $[e,1]$ with neutral element $e$, and let $\cl$ be a
closure operator on $L$.  Define $U_1$ by \eqref{eq:U1}.  Then the following
statements are equivalent:
\begin{enumerate}
\item $x\wedge e=\cl(x)\wedge e$ for every $x\in\Ie$;
\item conditions \textnormal{(C1)--(C4)} hold;
\item $U_1$ is a uninorm on $L$ with neutral element $e$;
\item $U_1$ is associative.
\end{enumerate}
Thus the closure equality is necessary and sufficient for this
bounded-lattice construction.
\end{corollary}

\begin{Proof}
The equivalence of (1) and (2) is \Cref{prop:lattice-reduction}, and
(2) implies (3) by \Cref{thm:U1}.  Clearly (3) implies (4).  It remains to
show that (4) implies (1).

Fix $x\in\Ie$ and write $c=c_e(x)$.  Since $x\le c$, one cannot have
$c\le e$.  Suppose that $c>e$.  Associativity at $(x,1,1)$ gives
$S(c,1)=c$, and hence $c=1$.  Associativity at $(x,x,1)$ then gives
\[
 x\wedge e=U_1\bigl(U_1(x,x),1\bigr)
 =U_1\bigl(x,U_1(x,1)\bigr)=1,
\]
which is impossible because $x\parallel e$.  Thus $c\in\Ie$.  Applying
associativity to $(1,x,c)$ gives
\[
 c\wedge e
 =U_1\bigl(U_1(1,x),c\bigr)
 =U_1\bigl(1,U_1(x,c)\bigr)
 =x\wedge e,
\]
where $x\le c$ was used in the last equality.  Finally,
\[
 c\wedge e
 =\bigl[\cl(x)\wedge(x\vee e)\bigr]\wedge e
 =\cl(x)\wedge e.
\]
Therefore $x\wedge e=\cl(x)\wedge e$, as required.
\end{Proof}

\begin{remark}
Ouyang and Zhang \cite[Theorem~4.1]{OuyangZhang2020} gave a sufficient
construction on a bounded lattice from a closure operator and a triangular
norm.  \c{C}ayl\i\ \cite[Theorems~3.1 and~3.11]{CayliKybernetika2023} studied
two further lattice formulas and obtained necessary and sufficient conditions.
Qi and Xiu \cite{QiXiu2025} instead start from a uninorm on a subinterval and
extend it by means of a closure operator and a triangular superconorm, while
Xiu and Si \cite{Xiu2025} consider closure operators together with triangular
conorms.  These papers are close in method, but their input data and defining
clauses are not identical in general.

The preceding corollary concerns only the operation in \eqref{eq:U1}.  For
this formula, it shows that the displayed closure equality is not merely
sufficient but also necessary.  It neither classifies all uninorms on a
bounded lattice nor subsumes all the constructions cited above.  On a proper
trellis, moreover, the ambient hypothesis in \Cref{thm:U1} cannot in general
be recovered from the lattice equality, as shown by the five-element example.
\end{remark}

\begin{corollary}\label{cor:lattice_case_U2}
Let $L$ be a bounded lattice, let $e\in L\setminus\{0,1\}$, let $T$ be a
triangular norm on $[0,e]$ with neutral element $e$, and let $\intr$ be an
interior operator on $L$.  Define $U_2$ by \eqref{eq:U2}.  Then the following
statements are equivalent:
\begin{enumerate}
\item $x\vee e=\intr(x)\vee e$ for every $x\in\Ie$;
\item conditions \textnormal{(I1)--(I4)} hold;
\item $U_2$ is a uninorm on $L$ with neutral element $e$;
\item $U_2$ is associative.
\end{enumerate}
\end{corollary}

\begin{Proof}
Apply \Cref{cor:lattice_case_U1} to the order-dual lattice.  Closure,
meet, triangular conorm, $c_e$, $m_e$, and \eqref{eq:U1} become,
respectively, interior, join, triangular norm, $d_e$, $n_e$, and
\eqref{eq:U2}.  This gives all four equivalences.
\end{Proof}

The constructions with an absorbing endpoint have the following corresponding
lattice criteria.

\begin{corollary}
\label{cor:lattice_case_U3_U4}
Let $L$ be a bounded lattice and let $e\in L\setminus\{0,1\}$.
\begin{enumerate}
\item Let $\cl$ be a closure operator, let $S$ be a triangular conorm on
$[e,1]$ with neutral element $e$, and define $U_3$ by \eqref{eq:U3}.
If $(e,1)\ne\varnothing$, then $U_3$ is a uninorm on $L$ with neutral
element $e$ if and only if $\cl(x)\wedge e=x\wedge e$ for every
$x\in\Ie$ and $S(x,y)\ne1$ for all $x,y\in(e,1)$.
If $(e,1)=\varnothing$, then $U_3$ is a uninorm for every closure
operator $\cl$.

\item Dually, let $\intr$ be an interior operator, let $T$ be a triangular
norm on $[0,e]$ with neutral element $e$, and define $U_4$ by
\eqref{eq:U4}.  If $(0,e)\ne\varnothing$, then $U_4$ is a uninorm on
$L$ with neutral element $e$ if and only if
$\intr(x)\vee e=x\vee e$ for every $x\in\Ie$ and $T(x,y)\ne0$ for all
$x,y\in(0,e)$.
If $(0,e)=\varnothing$, then $U_4$ is a uninorm for every interior
operator $\intr$.
\end{enumerate}
\end{corollary}

\begin{Proof}
Every element of $L$ satisfies (C1), and $m_e$ is jointly increasing.
Suppose first that $(e,1)\ne\varnothing$.  By
\Cref{prop:U3-relative-characterization}, the uninorm property of $U_3$
is equivalent to (C2), (C4), and condition
\eqref{eq:U3-relative-endpoint}.  Under (C1), condition (C3) follows from
(C2), while
\Cref{prop:lattice-reduction}(1) identifies the conjunction of (C2) and
(C4) with the displayed closure equality.  This proves the first
equivalence.  If
$(e,1)=\varnothing$, \Cref{prop:U3-relative-characterization} reduces
the criterion to joint increasingness of $m_e$, which is automatic on
a lattice.  The second assertion is the order dual, using
\Cref{cor:U4-relative-characterization} and
\Cref{prop:lattice-reduction}(2).
\end{Proof}

The following example treats the case in which the upper interval is
degenerate.

\begin{example}
\label{ex:U3-degenerate-upper-interval}
Let $L=\mathcal P(\{1,2,3\})$ be ordered by inclusion, with
$0=\varnothing$, $1=\{1,2,3\}$, and $e=\{1,2\}$.  Then
$(e,1)=\varnothing$ and
$\Ie=\bigl\{\{3\},\{1,3\},\{2,3\}\bigr\}$.
Define the closure operator by $\cl(A)=1$ for every $A\in L$.  For each
$x\in\Ie$ one has $c_e(x)=\cl(x)\cap(x\cup e)=1$,
so the image requirement in (C2) fails.  Nevertheless,
\Cref{cor:lattice_case_U3_U4} shows that $U_3$ is a uninorm.  Thus the
nonemptiness of $(e,1)$ is essential when the uninorm property of
$U_3$ is used to recover conditions \textnormal{(C1)--(C4)};
its order dual gives the corresponding boundary for $U_4$.
\end{example}

The two-sided constructions reduce to the following lattice criteria.

\begin{corollary}
\label{cor:lattice_case_U5_U6}
Let $L$ be a bounded lattice, let $e\in L\setminus\{0,1\}$, let $T$ be
a triangular norm on $[0,e]$, and let $S$ be a triangular conorm on
$[e,1]$, both with neutral element $e$.
\begin{enumerate}
\item Let $\cl$ be a closure operator and define $U_5$ by
\eqref{eq:U5}.  Then $U_5$ is a uninorm on $L$ with neutral element
$e$ if and only if $\cl(x)\parallel e$ for every $x\in\Ie$ and
condition \eqref{eq:zero-comparison-lower}
holds; explicitly,
\[
 p\le z,\quad p,q\in[0,e),\quad z\in\Ie
 \quad\Longrightarrow\quad T(p,q)=0.
\]
Equivalently, the first displayed condition may be replaced by
$\cl(x)\notin[e,1]$ for every $x\in\Ie$.  If $T$ has no zero divisors,
then condition \eqref{eq:zero-comparison-lower} may equivalently be replaced by
\eqref{eq:separation-lower}.

\item Dually, let $\intr$ be an interior operator and define $U_6$ by
\eqref{eq:U6}.  Then $U_6$ is a uninorm on $L$ with neutral element
$e$ if and only if $\intr(x)\parallel e$ for every $x\in\Ie$ and
condition \eqref{eq:one-comparison-upper}
holds; explicitly,
\[
 z\le p,\quad p,q\in(e,1],\quad z\in\Ie
 \quad\Longrightarrow\quad S(p,q)=1.
\]
Equivalently, the first displayed condition may be replaced by
$\intr(x)\notin[0,e]$ for every $x\in\Ie$.  If $S$ has no nontrivial
one-divisors, then condition \eqref{eq:one-comparison-upper} may equivalently be replaced by
\eqref{eq:separation-upper}.
\end{enumerate}
\end{corollary}

\begin{Proof}
Condition (C1) is automatic on $L$.  By
\Cref{prop:U5-relative-characterization}, $U_5$ is a uninorm if and
only if (C2) and \eqref{eq:zero-comparison-lower} hold.  Under (C1),
condition (C2) implies (C3); hence (C2) is equivalent here to
\textnormal{(C1)--(C3)}.  Part~(3) of \Cref{prop:lattice-reduction}
identifies these conditions with
$\cl(x)\parallel e$, equivalently $\cl(x)\notin[e,1]$, on $\Ie$.
This proves the first assertion.  The second assertion is its order
dual and follows from \Cref{cor:U6-relative-characterization}.  The final
equivalences follow from \Cref{cor:U5-no-zero-divisors} and its stated
order-dual counterpart.
\end{Proof}

\section{CONCLUSIONS}\label{sec:conclusion}

We have obtained necessary and sufficient conditions for six constructions on
bounded trellises.  Under the ambient transitivity condition,
the first closure-based construction is characterized by \textnormal{(C2)}
and \textnormal{(C4)}, while \textnormal{(C3)} follows from the other
assumptions; the interior-based result is its order dual.  For the variants
with an absorbing endpoint, the nondegenerate cases require
\eqref{eq:S_not_1} or \eqref{eq:T_not_0}, whereas in the degenerate cases the
closure or interior map no longer occurs in the effective formula.

For the constructions involving both a triangular norm and a triangular
conorm, the necessary and sufficient requirements are
\eqref{eq:zero-comparison-lower} and its order dual
\eqref{eq:one-comparison-upper}.  The stronger incomparability assumptions
\eqref{eq:separation-lower} and \eqref{eq:separation-upper} are equivalent to
these requirements when the triangular norm has no zero divisors or the
triangular conorm has no nontrivial one-divisors.  On bounded lattices, the
closure and interior parts of the criteria reduce to explicit identities or
image restrictions.

The examples also show the limitation of these reductions: a lattice
identity need not imply the ambient transitivity condition on a proper
trellis, and that condition need not be necessary for a particular operation
to be a uninorm.  Further work may determine weaker ambient hypotheses and
compare the resulting classes with other uninorm constructions on bounded
trellises.

\section*{DECLARATION OF COMPETING INTEREST}
The authors declare that they have no known competing financial interests or personal relationships that could have appeared to influence the work reported in this paper.

\section*{DATA AVAILABILITY}
No datasets were generated or analysed in this theoretical study.

\section*{ACKNOWLEDGEMENT}
\small
This work was supported by the National Natural Science Foundation of China
(Nos. 12671553 and 12271036).


\makecontacts

\end{document}